\documentclass[11pt]{article}

\usepackage{amssymb,amsmath,amsfonts,amssymb}
\hsize \textwidth
\advance \hsize by -\marginparwidth
\evensidemargin \oddsidemargin
\usepackage{amssymb}
\advance\hoffset by 5mm
\usepackage{graphics}

\def\@abssec#1{\vspace{.05in}\footnotesize \parindent .2in
{\bf #1. }\ignorespaces}

\newtheorem{theorem}{Theorem}[section]
\newtheorem{thm}[theorem]{Theorem}
\newtheorem{lemma}[theorem]{Lemma}

\newtheorem{corollary}[theorem]{Corollary}

\def \Rm {\mathbb R}

\newcommand{\cC}{{\mathcal C}}

\newcommand{\pdr}[2]{\dfrac{\partial{#1}}{\partial{#2}}}

\newcommand{\bx}{\mathbf x} 
 \newcommand{\vx}{\mathbf x}

\newcommand{\nwc}{\newcommand}
\nwc{\IE}{\mathbb{E}} 
\nwc{\IF}{\mathbb{F}} 
\nwc{\IP}{\mathbb{P}}
\nwc{\IG}{\mathbb{G}} 
\nwc{\IN}{\mathbb{N}} 
\nwc{\IQ}{\mathbb{Q}} 
\nwc{\IR}{\mathbb{R}} 
\nwc{\IT}{\mathbb{T}} 
\nwc{\IZ}{\mathbb{Z}} 
\nwc{\cH}{\mathcal H}

\title{Quenching of Reaction by Cellular Flows}
\author{Albert Fannjiang\thanks{Department of Mathematics,
University of California, Davis, CA 95616, USA; e-mail:
cafannjiang@ucdavis.edu} \and Alexander Kiselev\thanks{Institute
for Advanced Study, Princeton, NJ 08540 and Department of
Mathematics, University of Wisconsin, Madison, WI 53706, USA;
e-mail: kiselev@math.wisc.edu } \and Lenya Ryzhik
\thanks{Department of Mathematics, University of Chicago, Chicago,
IL 60637; e-mail: ryzhik@math.uchicago.edu}}

\begin{document}
\maketitle

\begin{abstract}
We consider a reaction-diffusion equation in a cellular flow. We
prove that in the strong flow regime there are two possible
scenario for the initial data that is compactly supported and the
size of the support is large enough. If the flow cells are large
compared to the reaction length scale, propagating fronts will
always form. For the small cell size, any finitely supported
initial data will be quenched by a sufficiently strong flow. We
estimate that the flow amplitude required to quench the initial
data of support $L_0$ is $A>CL_0^4\ln(L_0)$. The essence of the
problem is the question about the decay of the $L^\infty$ norm of
a solution to the advection-diffusion equation, and the relation
between this rate of decay and the properties of the Hamiltonian
system generated by the two-dimensional incompressible fluid flow.
\end{abstract}

\section{Introduction}

It has been  well understood since the classical work of G.I.
Taylor that the presence of a fluid flow may greatly increase the
mixing properties of diffusion. This phenomenon is known as ``eddy
diffusivity'' or ``enhanced diffusion''. The mathematical approach
to the problem is usually via the homogenization techniques that
concentrate on the long time-large scale behavior: see \cite{KM}
for a recent extensive review. This approach is appropriate when
there are no other time scales in the problem so that one may wait
as long as needed for the mixing effects to become prominent.

Recently there has been a lot of interest in the effect of flows
on the qualitative and quantitative behavior of solutions of
reaction-diffusion equations. Intuitively, there may be two
opposite effects of the additional mixing by the flow: on one
hand, it may increase the spreading rate of the chemical reaction
(the ``wind spreading the fire'' effect), or it may extinguish the
reaction (the ``try to light the campfire in a wind'' effect). The
first effect is related to the behavior of front-like solutions,
and has been extensively studied recently: traveling fronts have
been shown to exist in various flows
\cite{BH,BLL,BN,Volpert-2,X1,X2,X3}, and flows have been shown to
speed-up the front propagation due to the improved mixing
\cite{abp,B1,CKOR,KR,X4}, see \cite{B1,X3} for recent reviews of
the mathematical results in the area.  This problem has also
attracted a significant attention in the physical literature, we
mention
\cite{ACVV1,ACVV2,abp,KS,KRS,KBM,Kiss-Merkin1,Kiss-Merkin2} among
the recent papers and refer to \cite{Peters} as a general
reference. The present paper addresses the second phenomenon
mentioned above: the possibility of flame extinction by a flow.
The basic idea is that if the reaction process may occur only at
the temperatures $T$ above a critical threshold $\theta_0,$ then
mixing by a strong flow coupled to diffusion may drop the
temperature everywhere below $\theta_0$ and hence extinguish the
flame. However, unlike the usual linear advection-diffusion
homogenization problems, one may not wait for this to happen
beyond the time $t_c$ it takes for the chemical reaction to occur
-- the mixing has to happen before this time.  The question we
address is: "Given a threshold $\theta_0$, a time $t_c$, and the
support $L_0$ of the initial data, can we find a flow amplitude
$A_0(L_0)$ so that if the flow amplitude $A>A_0(L_0)$ is
sufficiently large then $\sup_{\vx}T(t=t_c,\vx)\le\theta_0$?" This
problem has been first considered in \cite{CKR} for
unidirectional, or shear, flows that have open streamlines. Even
in this simple situation the answer is non-trivial: in order for
quenching to be possible, the profile $u(y)$ should not be
constant on intervals larger than a prescribed size. The answer
has been shown to be sharp in \cite{KZ}.

In this paper, we study quenching by a different class of flows
with a more complex structure: incompressible cellular flows.
These are flows such that the whole plane $\Rm^2$ is separated
into invariant regions bounded by the separatrices of the flow
that connect the flow saddle points.  Many types of instabilities
in fluids lead to cellular flows, making them ubiquitous in
nature. We only mention Rayleigh-B\'enard instability in heat
convection, Taylor vortices in Couette flow between rotating
cylinders or heat expansion driven Landau-Darrieus instability.
The fact that the cellular flows have closed streamlines make the
effect of advection more subtle.  An important role in the
possibility of quenching is played by a thin boundary layer which
forms along the separatrices of the flow. Our main results show
that the cellular flow is quenching if and only if the size of the
minimal invariant regions (the flow cells) is smaller than a
certain critical size of the order of laminar flame length scale.

The simplest mathematical model that describes a chemical reaction in
a fluid is a single equation for temperature $T$ of the form
\begin{eqnarray}
  \label{eq:1.1}
  &&T_t+Au(\vx)\cdot \nabla T=\Delta T+M f(T)\\
&& T(0,\vx)=T_0(\vx)\nonumber
\end{eqnarray}
 where the flow $u(\vx)$ is prescribed.
We are interested in the effect of a strong advection, and accordingly
have written the velocity as a product of an amplitude $A$ and a fixed flow
$u(\vx)$. In this paper we consider nonlinearity of the ignition
type, that is, we assume that
\begin{eqnarray}
\nonumber
 && \hbox{(i) $f(T)$ is Lipschitz continuous on $0\le T\le 1$},\\
&&  \hbox{(ii) }~f(1)=0,~~\exists\theta_0
\hbox{ such that $f(T)=0$ for $T\in [0,\theta_0]$,
$f(T)>0$ for
$T\in (\theta_0, 1)$ }\label{eq:2.1.2}\\
&&\hbox{(iii) }~~f(T)\le T.\nonumber
\end{eqnarray}
The threshold $\theta_0$ is called the ignition temperature.
The last condition in (\ref{eq:2.1.2}) is just a normalization. We
consider the reaction-diffusion equation (\ref{eq:1.1}) in a two-dimensional
strip
$D=\left\{x\in {\mathbb R} ,~ y\in [0,2\pi l]\right\}$
with the periodic boundary conditions at the vertical boundaries:
\[
T(x,y+2\pi l)=T(x,y).
\]
The initial data $T_0(\vx)=T(0,\vx)$ is assumed to satisfy $0\le
T_0(\vx)\le 1$. The maximum principle implies that then $0\le
T(t,\vx)\le 1$ for all $t\ge 0$. We will say that regions with
temperature close to one are "hot", and those with temperature close
to zero are "cold".


The problem of extinction and flame propagation in (\ref{eq:1.1}) with
the ignition type nonlinearity (\ref{eq:2.1.2}) was first studied by
Kanel \cite{Kanel} in one dimension and with no advection. Assume for
simplicity that the initial data are given by a characteristic
function: $T_0(\vx)=\chi_{[0,L]}(x).$ Kanel showed that, in the
absence of fluid motion, there exist two length scales $L_0<L_1$ such
that the flame becomes extinct for $L<L_0$, and propagates for
$L>L_1$.  More precisely, he has shown that there exist $L_0$ and
$L_1$ such that
\begin{eqnarray}
  \label{eq:2.2}
&&  T(t,\vx)\to 0~\hbox{as $t\to\infty$ uniformly in $D$ if $L<L_0$}\\
&&T(t,\vx)\to 1~\hbox{as $t\to\infty$ for all $(x,y)\in D$ if
$L>L_1$}.
\nonumber
\end{eqnarray}
In the absence of advection, the flame extinction is achieved by
diffusion alone, given that the support of initial data is small
compared to the scale of the laminar front width $l_c =M^{-1/2}.$
However, in many applications quenching is a result of strong wind,
intense fluid motion and operates on larger scales. Kanel's result was
extended to non-zero advection by shear flows by Roquejoffre
\cite{Roq-2} who has shown that (\ref{eq:2.2}) holds also for $u\ne 0$
with $L_0$ and $L_1$ depending, in particular, on $A$ and $u(y)$ in an
uncontrolled way.

As we have mentioned, the question of the dependence of the strength of
advection $A$ which is necessary for quenching the initial data of a given
size $L_0$ has been recently studied in \cite{CKR} and \cite{KZ}
in the case of a unidirectional (shear) flow $(Au(y),0)$.
Following \cite{CKR}, we call the flow $u(x,y)$ quenching if for
every $L_0$ there exists $A_0(L_0)$ such that the solution of
(\ref{eq:1.1}) with the initial data of size $L_0$ and advection
strength $A>A_0(L_0)$ quenches.
It turns out \cite{CKR,KZ} that the shear flow $u(y)$ is quenching
if and only if $u(y)$ does not have a plateau of size larger than
a certain critical threshold (comparable with the length scale
$M^{-1/2}$ which characterizes the width of a laminar flame). The
intuition behind this result is that shear flows are very
effective in stretching the front and exposing the hot initial
data to cool-off effects of diffusion unless there is a
long,
flat part in their profile, where this phenomenon
is obviously not present.

\begin{figure}
\begin{center}
\scalebox{0.73}{\includegraphics{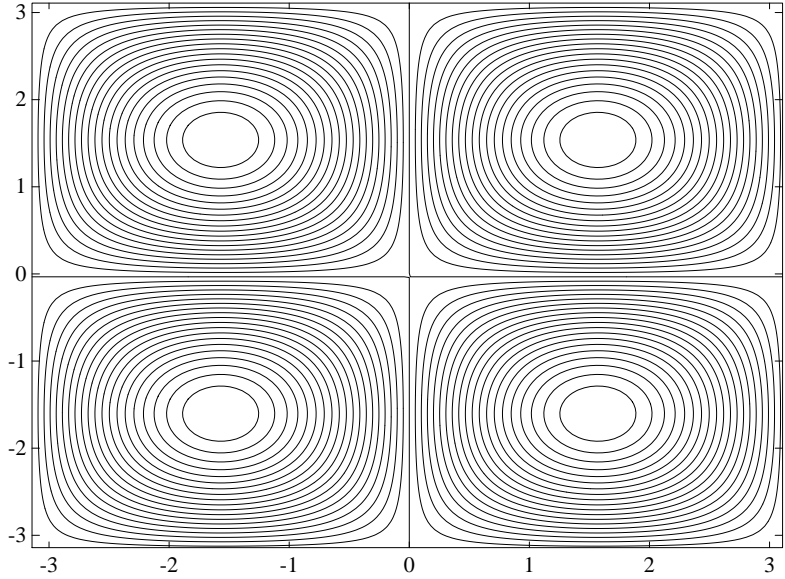}} \caption{Cellular
flow} \label{struct}
\end{center}
\end{figure}


Here we consider (\ref{eq:1.1}) for the domain $D=\IR\times [-\pi
l, \pi l]$ with the $2\pi l$-periodic boundary conditions in $y$
and decay conditions in $x$:
\begin{equation}\label{BC}
T(t,x,y)=T(t,x,y+2\pi l),~~~T(t,x,y)\to 0~~\hbox{as
$x\to\pm\infty$.}
\end{equation}
We restrict ourselves to a particular example of a cellular flow
$u(x,y)$ that has the form $u_l(x,y)=\nabla^\perp h_l(x,y)$, where
$\nabla^\perp h=(h_y,\!-h_x)$. Here $l$ defines the size of a flow
cell, and we take the stream function $h_l$ to be
\begin{equation}\label{intro-hl}
h_l(x,y)= l\sin\frac{x}{l} \sin
\frac{y}{l}
\end{equation}
whose streamline structure are shown in Fig. 1.

We will usually omit the index $l$ in notation for $u_l$,
$h_l$ for the fluid flow; it will be clear from the context what
the scaling is. The initial data $T_0(x,y)$ is non-negative, and
bounded above by one: $0\le T_0(x,y)\le 1$.

The first theorem shows that cellular flows with large cells do not have
the quenching property.
\begin{thm}\label{thm1} Assume that $T_0(x,y)=1$ for $(x,y)\in
[0,\pi l]\times[0,\pi l]$. There exists a critical cell size $l_0
\sim M^{-1/2}$ so that if $l\ge l_0$, then under the above
assumptions on the advection $u$, we have $T(t,x,y)\to 1$ as $t\to
+\infty$, uniformly on compact sets, for all $A\in{\mathbb R}$.
\end{thm}
The notation $l_0 \sim M^{-1/2}$ means that $C_1 M^{-1/2} \leq l_0
\leq C_2 M^{-1/2}$ with $C_{1,2}$ some positive universal constants.
The proof of Theorem \ref{thm1} is simple and is based on a
construction of an explicit time-independent sub-solution.

Next, we show that if the cell size $l$ is sufficiently small then a
sufficiently strong flow will quench a flame. More precisely, we have
the following result.
\begin{thm}\label{thm2} Assume that $T_0(x,y)=0$ outside an interval
$-L_0\le x\le L_0$ and $0\le T_0(x,y)\le 1$ for all $(x,y)\in D$.
There exists a critical cell size $l_0 \sim M^{-1/2}$ so that if
$l\le l_0$, then there exists $A_0(L_0)$ such that we have
$T(t,x,y)\to 0$ as $t\to +\infty$, uniformly in $D$, for all $A\ge
A_0(L_0)$. For large $L_0,$ we have $A_0(L_0) \leq
C(l)L_0^4\ln(L_0).$
\end{thm}



 A formal argument
based on the homogenization theory predicts that $A_0\sim L_0^4$
without the factor of $\ln L_0$ -- this follows from the effective
diffusivity scaling $\kappa_*\sim\sqrt{A}$ that was first shown
formally in \cite{Childress} and later proved in
\cite{FP,Koralov,NPR}. The same scaling may be obtained from the
formal predictions $V_A\sim A^{1/4}$ for the front speed $V_A$ in
a cellular flow \cite{ACVV1,ACVV2,abp,VCKRR} -- this implies that
the front width is of the order $A^{1/4}$.  Hence one might expect
that initial data with the support less than the front width
$L_0<A^{1/4}$ to be quenched. Therefore, the rigorously proved
bound of Theorem~\ref{thm2} is likely to be sharp up to a
logarithmic factor.



The particular choice of the stream function (\ref{intro-hl}) is
not important for the proof but it does simplify some of the
estimates -- it is straightforward to generalize our result to
other cellular flows. However, the proof of Theorem \ref{thm2}
does use the periodicity of the flow in an essential way. We
believe that the quenching property should hold for sufficiently
regular non-periodic flows with small cells as well. We have
preliminary results in this direction using different techniques;
however, these results give much weaker upper bound for
$A_0(L_0),$ and will appear elsewhere.


At the heart of the proof of Theorem \ref{thm2} is the question
about the rate of decay of $L^\infty$ norm (and its dependence on
$A$) of the solution to a passive advection-diffusion equation.
Thus, our main object of study is a natural question about the
effects of a combination of two fundamental and separately
well-understood processes: advection by a fixed incompressible
flow and diffusion. Yet their interaction is well known to produce
subtle phenomena. The issues we study are directly related to the
work of Freidlin and Wentzell \cite{Freidlin,FW-0,FW} on the
random perturbations of the Hamiltonian systems. They show that in
the limit of large $A$ the process converges to a diffusion on the
Reeb graph of the background Hamiltonian. The relation to that
work is very natural as any incompressible flow in two dimensions
is a Hamiltonian system (the stream function is the Hamiltonian).
However, the cellular flow that we consider does not satisfy the
assumptions of the Freidlin-Wentzell theory, which requires growth
of the Hamiltonian at infinity and does not allow existence of
hetero-clinic orbits. Nevertheless, one may restate the small cell
assumption in Theorem \ref{thm2} as a requirement that the Reeb
graph of the Hamiltonian has a sufficiently small diameter in its
natural metric. The proof of Theorem \ref{thm2} is based on two
observations: first, temperature becomes approximately constant on
the whole skeleton of separatrices, as the skeleton is just one
point on the Reeb graph. Because of that, solution inside each
cell may be split into two parts. One solves the initial value
problem with zero data on the boundary. Another solves a nearly
identical boundary value problem in each cell. The first part
decays because the cell is small -- this is where we use the size
restriction. The second one is nearly identical on all cells,
hence it has to be small in order not to violate the preservation
of the $L^1$-norm. The technical part of the proof is in making
this informal scenario rigorous.

The paper is organized as follows: in Section~\ref{sec:absence} we
prove Theorem~\ref{thm1} by constructing an appropriate sub-solution
and using certain PDE estimates to prove convergence of solution to
unity. In Section~\ref{sec:quench} we give the proof of
Theorem~\ref{thm2}, which is more involved. It uses uniform in $A$
estimates on the evolution of advection-diffusion equation, a boundary
layer argument, and probabilistic estimates for an auxiliary cell
heating problem.

\noindent {\bf Acknowledgment.}
The research of AF is supported in part by The
Centennial Fellowship from American Mathematical
Society
and U.S. National Science Foundation (NSF) grant DMS-0306659. AK has been
supported in part by NSF grants DMS-0321952 and
DMS-0314129. LR has been supported in part by NSF
grant DMS-0203537 and ONR grant N00014-02-1-0089. Both
AK and LR acknowledge support by Alfred P. Sloan
fellowships.

\section{Absence of quenching by large cells}\label{sec:absence}

We prove in this section Theorem \ref{thm1}, that is, we show that
cellular flows with sufficiently large cells do not have the quenching
property.
The proof consists of two steps. First, we construct a
time-independent sub-solution $\Phi(\vx)$ to (\ref{eq:1.1}) in a
cell ${\mathcal C}_1=[0,\pi l]\times[0,\pi l]$: the function
$\Phi$ satisfies
\begin{equation}\label{sub-cond-0}
{A}u_l\left(\vx\right)\cdot\nabla\Phi\le\Delta\Phi+ Mf(\Phi)
\end{equation}
and is $2\pi l$-periodic in $y$.
It is also positive on an open set inside ${\mathcal C}_1$ and
negative on $\partial {\mathcal C}_1$. We normalize $\Phi$ so that $\Phi\le 1$.
As $T_0=1$ on ${\mathcal C}_1$ by assumption, we have
$T_0(\vx)\ge\Phi(\vx)$. Then the maximum principle implies that $T(t,\vx)\ge \Phi(\vx)$ for
all $t\ge 0$ and $\vx\in {\mathcal C}_1$. It follows that
$T(t,\vx)$ does not vanish as $t\to +\infty$. In the second step
we show that actually $T(t,\vx)\to 1$. We begin with the construction of
the sub-solution $\Phi(\vx)$. First, we rescale equation
(\ref{eq:1.1}) by $\vx \to l\vx$ so that a
sub-solution in the rescaled coordinates should
satisfy
\begin{equation}\label{sub-cond}
\frac{A}{l}u(\vx)\cdot\nabla\Phi\le\frac{1}{l^2}\Delta\Phi+ M
f(\Phi).
\end{equation}
\begin{lemma}\label{lemma1}
  If $l$ is sufficiently large then there exists a $C^1$ function
  $\Phi(\vx)$ that is constant on the streamlines of the flow $u(\vx)$
  and satisfies (\ref{sub-cond}) for $\vx\in
  {\mathcal C}_1=[-\pi,\pi]^2$.  Moreover,
$\Phi(\vx)\le 1$ for all
  $\vx\in {\mathcal C}_1$, and $\Phi(\vx)<0$ for $\vx\in\partial {\mathcal C}_1$.
\end{lemma}
{\bf Proof of Lemma \ref{lemma1}.} We may choose two numbers $\theta_1$ and
$\theta_2$ so that $\theta_0<\theta_1<\theta_2<1$ and such that the
straight line that connects the point $(\theta_1,0)$ to the point
$(\theta_2,f(\theta_2))$ lies below the graph of $f(T)$. More precisely, that
means that the function
\begin{equation}\label{def-g}
g(T)=\left\{\begin{matrix} 0, & T\le\theta_1\cr
                          \alpha(T-\theta_1), & \theta_1\le T\le \theta_2\cr
                          f(T),& T\ge\theta_2\cr\end{matrix}\right.
\end{equation}
with $\alpha=f(\theta_2)/(\theta_2-\theta_1)$, satisfies $g(T)\le
f(T)$.  Such modification of $f(T)$ for the construction of sub-solutions
was first used in \cite{Kanel}, and then in \cite{CKR}. A function
$\Phi(h(\vx))$ satisfies (\ref{sub-cond}) if
\[
\frac{1}{l^2}\left(|\nabla h|^2\frac{d^2\Phi}{d h^2}+ \Delta
h\frac{d\Phi}{dh}\right)+ M g(\Phi)\ge 0,
\]
where $h(x,y)=\sin x\sin y$ is the stream function. Note that the
advection term vanishes identically for such functions. We   use the fact that $\Delta h=-2h$
to obtain
\begin{equation}\label{indicate}
|\nabla h|^2\frac{d^2\Phi}{d h^2}-2h\frac{d\Phi}{dh}+
\left(\frac{l}{l_c}\right)^2g(\Phi)\ge 0.
\end{equation}
Here $l_c = M^{-1/2}$ is the laminar front width, the length scale
associated to the chemical reaction strength. Relation
(\ref{indicate}) indicates that the ratio $l/l_c$ has to be
sufficiently large for a sub-solution to exist.
Note that
\[
2h(1-h)\le |\nabla h|^2\le 2(1-h^2).
\]
Indeed, we have
\begin{eqnarray*}
&&|\nabla h(x,y)|^2=
\sin^2x\cos^2y+\cos^2x\sin^2y=\sin^2x+\sin^2y-2\sin^2x\sin^2y\\
&&\ge 2\sin x\sin y(1-\sin x\sin y)=2h(x,y)(1-h(x,y))
\end{eqnarray*}
and
\begin{eqnarray*}
|\nabla h(x,y)|^2
\le 2-2\sin^2x\sin^2y=2(1-h^2(x,y)).
\end{eqnarray*}
Therefore it suffices to construct an increasing function
$\Phi(h)$ satisfying
\[
\frac{d^2\Phi}{dh^2}-\frac{2h}{2h(1-h)}\frac{d\Phi}{dh}+
\frac{1}{2(1-h^2)}\left(\frac{l}{l_c}\right)^2g(\Phi)\ge 0,
\]
which would in turn follow from
\[
\frac{d^2\Phi}{dh^2}-\frac{1}{1-h}\frac{d\Phi}{dh}+
\frac{1}{2}\left(\frac{l}{l_c}\right)^2g(\Phi)= 0.
\]
Make a change of variables $\displaystyle
R=\frac{l}{l_c\sqrt{2}}(1-h)$ so that the above becomes
\begin{equation}\label{bessel-eq}
\frac{d^2\Phi}{dR^2}+\frac{1}{R}\frac{d\Phi}{dR}+g(\Phi)=0.
\end{equation}
The center of the cell corresponds now to $R=0$, while the
boundary $h=0$ becomes $R=\displaystyle\frac{l}{l_c\sqrt{2}}$.  We
impose the following ``initial data'' for (\ref{bessel-eq}):
\[
\Phi(0)=\theta_2,~~\frac{d\Phi(0)}{dR}=0.
\]
The explicit form (\ref{def-g}) of the function $g(\Phi)$
implies that the solution $\Phi(R)$ is given explicitly by
\begin{eqnarray}
  \label{eq:3.7}
&& \Phi(R)=\theta_1+(\theta_2-\theta_1)
J_0\left(R\sqrt{\alpha}\right),~~
\hbox{ for $\displaystyle R\le R_1=\frac{\xi_1}{\sqrt{\alpha}}$}
\end{eqnarray}
with $\alpha$ as in (\ref{def-g}).
Here $J_0(\xi)$ is the Bessel function of order zero, and $\xi_1$ is
its first zero.

Furthermore, we have
\begin{eqnarray}
  \label{eq:3.8}
  \Phi(R)=B\ln\frac{R_2}{R},~~ \hbox{for $R_1\le R$}.
\end{eqnarray}
The constants $B$ and $R_2$ are determined by matching the functions
(\ref{eq:3.7}) and (\ref{eq:3.8}), and their derivatives at
$R=R_1$. Then we get
\begin{eqnarray*}
B=(\theta_2-\theta_1)\xi_1|J_0'(\xi_1)|,~~
  R_2=\xi_1\sqrt{\frac{\theta_2-\theta_1}{f(\theta_2)}}
\exp\left[\frac{\theta_1}{(\theta_2-\theta_1)\xi_1 |J_0'(\xi_1)|}\right].
\end{eqnarray*}
Observe that the function $\Phi(R)$ constructed above is negative on the boundary
of the cell only provided that $R_2<\displaystyle\frac{l}{l_c\sqrt{2}}$, which means
that the cell size
\[
l\ge l_c\sqrt{2}\xi_1\sqrt{\frac{\theta_2-\theta_1}{f(\theta_2)}}
\exp\left[\frac{\theta_1}{(\theta_2-\theta_1)\xi_1 |J_0'(\xi_1)|}\right]
\]
has to be sufficiently large for this construction to be applicable.
This proves Lemma \ref{lemma1}. $\Box$

In order to finish the proof of Theorem \ref{thm1} we have to show
that $T(t,\vx)\to 1$ as $t\to +\infty$ provided that $T_0(\vx)=1$ on
a cell. For such initial data we have
\[
T_0(\vx)\ge \Phi_0(\vx)=\max\left\{\Phi(h(\vx)),0\right\},
\]
where the function $\Phi(\vx)$ is the sub-solution constructed in
Lemma \ref{lemma1}.  It follows from the parabolic maximum principle that
then $T(t,\vx)\ge \Phi_0(\vx)$ for
all $t\ge 0$. Furthermore, we have $T(t,\vx)\ge \Psi(t,\vx)$,
where the function $\Psi(t,\vx)$ satisfies (\ref{eq:1.1}) with the
initial data $\Phi_0(\vx)$.  Note that $\Psi(t,\vx)\ge\Phi_0(\vx)$
for all $t\ge 0$. The maximum principle applied to the finite differences
$\Psi_h(t,x)=\Psi(t+h,\vx)-\Psi(t,\vx)$
implies that $\Psi(t,\vx)$
is a point-wise increasing function of time that is bounded above
by one. Therefore the point-wise limit
$\overline\Psi(\vx)=\lim_{t\to+\infty}\Psi(t,\vx)$ exists,
moreover, $\overline\Psi(\vx)\ge\Phi_0(\vx)$, and
$\overline\Psi(\vx)$ satisfies the stationary problem
\begin{eqnarray}\label{Psibar-eq}
&&A u(\vx)\cdot\nabla \overline\Psi=\Delta\overline\Psi +M
f(\overline\Psi).
\end{eqnarray}
with the $2\pi l$-periodic boundary conditions in $y$
\[
\overline\Psi(t,x,0)= \overline\Psi(t,x,2\pi l)
\]
(here $2N$ is the number of cells in $y$ direction).
We have the following lemma.
\begin{lemma}\label{lemma2} The function $\overline\Psi(x,y)$ satisfies the
following bound:
\begin{equation}\label{eq-lemma2}
\int_{D}
|\nabla \overline\Psi(\vx)|^2 \,d\vx+\int_D f(\overline\Psi(\vx))\,d\vx<+\infty,
\end{equation}
where $D={\mathbb R}_x\times[-\pi l,\pi l]_y$.
\end{lemma}
{\bf Proof.} The function $\Psi(t,\vx)$ satisfies an a priori bound
\[
\frac{1}{\tau}\int_0^\tau \left(\int_{D} f(\Psi(t,\vx))\,d \vx
\right)dt\le C_0,~~~ \tau\ge CM^{-1},
\]
that may be easily proved as in \cite{CKOR}. Here the constant $C_0$ may depend
on the flow amplitude $A$. Therefore, there
exists a sequence of times $t_n\to +\infty$ so that
\[
\int_{D} f(\Psi(t_n,\vx))\,d\vx\le C_0.
\]
This implies that
\[
\int_{D} f(\overline\Psi(\vx))\,d\vx\le C_0,
\]
and it remains to obtain the bound on $\|\nabla \overline\Psi\|_{L^2}$
in (\ref{eq-lemma2}). We multiply (\ref{Psibar-eq}) by $\overline\Psi$
and integrate in $x$ between $-X+\zeta$ and $X+\zeta$ with $X$ large and
$\zeta\in[0,l_c]$, and in $y\in{\mathbb R}$.  We get
\begin{eqnarray*}
&& \frac{A}{2}\int_{-\pi l}^{\pi l}\left[u_1(X+\zeta,y)
|\overline\Psi(X+\zeta,y)|^2-
u_1(-X+\zeta,y)|\overline\Psi(-X+\zeta,y)|^2\right]dy\\ &&=
\int_{-\pi l}^{\pi l}\left[
\overline\Psi(X+\zeta,y)\overline\Psi_x(X+\zeta,y)-
\overline\Psi(-X+\zeta,y)\overline\Psi_x(-X+\zeta,y)\right]dy\\
&&+M \int_{-\pi l}^{\pi l}dy\int_{-X+\zeta}^{X+\zeta}
\overline\Psi(x,y)f(\overline\Psi(x,y))dx- \int_{-\pi l}^{\pi
l}dy\int_{-X+\zeta}^{X+\zeta} |\nabla\overline\Psi(x,y)|^2dx
\end{eqnarray*}
and average this equation in $\zeta\in[0,l_c]$. This provides the bound
\begin{eqnarray*}
\frac{1}{l_c}\int_0^{l_c}d\zeta\int_{-\pi l}^{\pi l}
dy\int_{-X+\zeta}^{X+\zeta} |\nabla\overline\Psi(x,y)|^2dx\le \pi
Al\|u\|_{\infty}+ \frac{2\pi l}{l_c}+ M
\int_{D}f(\overline\Psi(x,y))dxdy
\end{eqnarray*}
Taking $X$ to infinity, we obtain
\begin{eqnarray*}
\int_D |\nabla\overline\Psi(x,y)|^2dx\le \pi Al\|u\|_{\infty}+
\frac{2\pi l}{l_c}+ M \int_{D}f(\overline\Psi(x,y))dxdy
\end{eqnarray*}
and the bound on $\nabla\overline\Psi$ in (\ref{eq-lemma2})
follows.
$\Box$

Now we are ready to complete the proof of Theorem~\ref{thm1}. \\
{\bf Proof.} Lemma \ref{lemma2} implies that there exist two sequences of points
$x_n\to -\infty$ and $z_n\to+\infty$ so that
\begin{equation}\label{def-xnzn}
\int_{-\pi l}^{\pi l}\left(|\nabla\overline\Psi(x_n,y)|^2+
|\nabla\overline\Psi(z_n,y)|^2\right)dy\to 0~~\hbox{as $n\to +\infty$.}
\end{equation}
We integrate (\ref{Psibar-eq}) in $y$ and in $x$ between $x_n$ and
$z_n$ to obtain
\begin{eqnarray}\label{xn-zn}
&&{A}\int_{-\pi l}^{\pi
l}\left[u_1(z_n,y)\overline\Psi(z_n,y)-
u_1(x_n,y)\overline\Psi(x_n,y)\right]dy\\
&&=\int_{-\pi l}^{\pi l}\left[
\overline\Psi_x(z_n,y)-\overline\Psi_x(x_n,y)\right]dy +M
\int_{-\pi l}^{\pi l}dy\int_{x_n}^{z_n}
f(\overline\Psi(x,y))dx.\nonumber
\end{eqnarray}
We pass to the limit $n\to\infty$ in (\ref{xn-zn}). Observe that
\[
\left|\int_{-\pi l}^{\pi l}
\overline\Psi_x(z_n,y)dy\right|\le
\sqrt{\pi l}\left(\int_{-\pi l}^{\pi l}
|\overline\Psi_x(z_n,y)|^2dy\right)^{1/2}
\to 0 ~~\hbox{as $n\to +\infty$}
\]
as follows from (\ref{def-xnzn}),
and similarly
\[
\left|\int_{-\pi l}^{\pi l}
\overline\Psi_x(x_n,y)dy\right|\to 0~~\hbox{as $n\to +\infty$.}
\]
Furthermore, since
\[
\int^{\pi l}_{-\pi l} u_1(x,y)dy=0
\]
for all $x\in{\mathbb R}$, (\ref{def-xnzn}) and the Cauchy-Schwartz
inequality imply that
\[
\int_{-\pi l}^{\pi l}u_1(z_n,y)\overline\Psi(z_n,y)dy=
\int_{-\pi l}^{\pi l}u_1(z_n,y)\left(
\int_0^y\overline\Psi_\xi(z_n,\xi)d\xi\right)d y\to 0,
~~\hbox{as $n\to +\infty$.}
\]

Therefore in the limit $n\to +\infty$ equation (\ref{xn-zn}) becomes
\[
\int_{D}f(\overline\Psi(x,y))dxdy=0
\]
and hence $f(\overline\Psi(x,y))=0$ for all $(x,y)\in D$. However,
since $\overline\Psi(x,y)\ge\Phi_0(x,y)$, while, on the other hand,
$\max\Phi_0(x,y)=\theta_2>\theta_0$, and the function $
\overline\Psi(x,y)$ is continuous, we conclude that
$\overline\Psi(x,y)\equiv 1$. This finishes the proof of Theorem
\ref{thm1}. $\Box$

\section{Quenching by small cells}\label{sec:quench}

In this section we show that the cellular flow with small cells is
quenching.
The proof proceeds in several steps.
First, we reduce the problem to a linear advection-diffusion
equation. Indeed, as $f(T)\le T$ we have the following upper
bound for $T:$
\begin{equation}\label{small-Tphi}
T(t,\vx)\le e^{M t}\phi(t,\vx).
\end{equation}
The function $\phi(t,\vx)$ satisfies the advection-diffusion equation
\begin{equation}\label{small-linear}
\pdr{\phi}{t}+Au\cdot\nabla \phi= \Delta\phi
\end{equation}
with the same initial data $\phi(0,\vx)=T_0(\vx)$ and the
$2\pi l$-periodic boundary conditions in $y$:
$\phi(t,x,y)=\phi(t,x,y+2\pi l)$. Note that if at some time $t_0>0$ we
have $\phi(t_0,\vx)\le\theta_0$ everywhere, then the
maximum principle implies that
$\phi(t,\vx)\le\theta_0$ and $T$ satisfies the linear
equation (\ref{small-linear}) for all $t\ge t_0$. Then
the conclusion of Theorem \ref{thm2} follows. Hence,
the upper bound (\ref{small-Tphi}) implies that it
suffices to show
$\phi(t=M^{-1},x,y)\le\theta_0'=\theta_0e^{-1}$ and this is what we will do.

Heuristically, the proof relies on the observation that solution
of (\ref{small-linear}) should generally become constant along the
streamlines of the flow if its amplitude is large.  Moreover, the
value of the solution on the streamlines $h=h_0$ very near the
boundary in two neighboring cells have to be close (as follows
from a simple $L^2$-bound on $\nabla\phi$ appearing in
Lemma~\ref{streamconst} below). However, that means that solution
should have, roughly speaking, the same profile in each cell. This
is incompatible with the preservation of the $L^1$-norm of $\phi$
unless this function is very small in each of the cells which
means that solution has to be less than $\theta_0$ everywhere. The
proof follows this heuristic outline -- the technical difficulty
is that we are able to control the uniformity of the solution
along the streamlines only in a space-time averaged sense.
Additional ingredients are required to obtain the point-wise
control.

\subsection{The  Nash inequality lemma}

We will need throughout the proof an $L^1-L^\infty$ decay estimate
for solutions of the linear diffusion-advection
\begin{equation}\label{small-adv-diff}
\pdr{T}{t}+v\cdot\nabla T=\Delta T
\end{equation}
that is independent of the advection strength.
Equation (\ref{small-adv-diff}) is considered in the infinite strip $D
=\IR\times [-\pi l, \pi l]$
with the $2\pi l$-periodic boundary conditions in $y$
direction.
\begin{lemma}\label{lemma-l1infty}
There exists a constant $C>0$ so that the solution of
\begin{eqnarray}\label{v-eq}
&&\pdr{\psi}{t}+v\cdot\nabla\psi=\Delta\psi\\
&&\psi(0,\vx)=\psi_0(\vx)\ge 0,~~~\vx\in{\mathbb R}^2\nonumber
\end{eqnarray}
with the $2\pi l$-periodic boundary condition in $y$ and a flow
$v$ that is
$2\pi l$-periodic, sufficiently regular and divergence-free:
$\nabla\cdot v=0$, satisfies
\begin{equation}\label{l1-linfty}
\|\psi(t)\|_{L^\infty(D)}\le
{C}n^2(t)\|\psi_0\|_{L^1(D)},
\end{equation}
where $D={\mathbb R}_x\times[-\pi l,\pi  l]_y$.
Here $n(t)$ is the unique solution of
\begin{equation}\label{n-eq}
\frac{4n^4(t)}{1+4n^3(t)l^3}=\frac{C_1}{ l^2t,}
\end{equation}
and the constants $C,C_1$ do not depend on $v.$
\end{lemma}
\it Remark. \rm Note that (\ref{n-eq}) implies that for $t \geq
l^2,$ we have
$n(t)^2 \sim \frac{C}{l \sqrt{ t}}.$ Hence, solution decays as
the solution of the one-dimensional problem after the  time it takes the
diffusion to feel the boundary.

{\bf Proof.}
We multiply (\ref{v-eq}) by $\psi$ and integrate over the domain $D$ to obtain
\begin{equation}\label{2.9}
\frac{1}{2}\frac{d}{dt}\|\psi\|_2^2=-\|\nabla\psi\|_2^2.
\end{equation}
Here and below in the proof of this Lemma, $\|\cdot\|_p$ denotes the
norm in $L^p(D)$.

We now prove the following version of the Nash
inequality \cite{Nash} for a strip of width $l$:
\begin{equation}\label{nash-strip}
\|\nabla\psi\|_2^2\ge
C\frac{l^2\|\psi\|_2^6}{\|\psi\|_1^4+l^3\|\psi\|_1\|\psi\|_2^3}
\end{equation}
The proof of (\ref{nash-strip}) is similar to that of the usual Nash
inequality. We represent $\psi$ in terms of its Fourier series-integral:
\[
\psi(x,y)=\sum_{n\in{\mathbb Z}}\int_{\mathbb R} e^{ iny/l+ikx}
\hat \psi_n(k)\frac{dk}{(2\pi)^2},
\]
where
\[
\hat \psi_n(k)=\frac{1}{ l}\int_D e^{-ikx- iny/l}\psi(x,y)dxdy.
\]
Therefore we have $|\hat \psi_n(k)|\le \displaystyle\frac{1}{l}\|\psi\|_{L^1}$.
The Plancherel formula becomes
\begin{eqnarray*}
&&\int_D|\psi(x,y)|^2dxdy=\sum_{n,m\in{\mathbb Z}}
\int_{{\mathbb R}^2} e^{ iny/l-imy/l
+ikx-ipx}\hat
\psi_n(k)\overline{\hat{\psi}}_m(p)\frac{dkdpdxdy}{(2\pi)^4}\\
&&~~~~~~~~~~~~~~~~~~~~~~~= l\sum_{n\in{\mathbb Z}}\int|\hat
\psi_n(k)|^2\frac{dk}{(2\pi)^2}
\end{eqnarray*}
and similarly
\begin{eqnarray*}
&&\int_D|\nabla\psi(x,y)|^2dxdy=
l\sum_{n\in{\mathbb Z}}\int_{\mathbb
R}\left(k^2+\frac{n^2}{l^2}\right) |\hat
\psi_n(k)|^2\frac{dk}{(2\pi)^2}.
\end{eqnarray*}
Let $\rho>0$ be a positive number to be chosen later. Then using the Plancherel
formula we may write
\[
\|\psi\|_2^2=I+II,
\]
where
\[
I=l\sum_{|n|\le \rho l}\int_{|k|\le \rho}|\hat \psi_n(k)|^2\frac{dk}{2\pi}\le
\frac{Cl\rho([l\rho]+1)}{l^2}\|\psi\|_1^2\le
\frac{C\rho(l\rho+1)}{l}\|\psi\|_1^2.
\]
The rest may be bounded by
\[
II\le \frac{l}{\rho^2} \sum_{n\in{\mathbb Z}}\int_{k\in{\mathbb
R}} \left(k^2+\frac{n^2}{l^2}\right)
|\hat{\psi}_n(k)|^2\frac{dk}{(2\pi)^2}\le
\frac{C}{\rho^2}\|\nabla\psi\|_2^2.
\]
Therefore we have for all $\rho>0$:
\[
\|\psi\|_2^2\le  \frac{C\rho(l\rho+1)}{l}\|\psi\|_1^2+
\frac{C}{\rho^2}\|\nabla\psi\|_2^2.
\]
We choose $\rho$ so that
\[
\rho^3=\frac{l\|\nabla\psi\|_2^2}{\|\psi\|_1^2}
\]
and obtain
\begin{eqnarray*}
&&\|\psi\|_2^2\le \frac{C\|\nabla\psi\|_2^{2/3}}{l^{2/3}\|\psi\|_1^{2/3}}
\left(\frac{l^{4/3}\|\nabla\psi\|_2^{2/3}}{\|\psi\|_1^{2/3}}+1\right)
\|\psi\|_1^2+\frac{C\|\nabla\psi\|_2^{2}\|\psi\|_1^{4/3}}{l^{2/3}
\|\nabla\psi\|_2^{4/3}}\\
&&~~~~~=
\frac{2C}{l^{2/3}}\|\psi\|_1^{4/3}\|\nabla\psi\|_2^{2/3}+
Cl^{2/3}\|\nabla\psi\|_{2}^{4/3}\|\psi\|_1^{2/3}.
\end{eqnarray*}
This is a quadratic inequality $ax^2+bx-c\ge 0$ with
$x=\|\nabla\psi\|_2^{2/3}$, $a=l^{2/3}\|\psi\|_1^{2/3}$,
$b=\displaystyle\frac{2}{l^{2/3}}\|\psi\|_1^{4/3}$, and $c=\|\psi\|_2^2/C$
and hence
\[
x\ge\frac{-b+\sqrt{b^2+4ac}}{2a}=\frac{2c}{b+\sqrt{b^2+4ac}}\ge
\frac{c}{\sqrt{b^2+4ac}}.
\]
This implies that
\[
\|\nabla\psi\|_2^{2/3}\ge
C{\|\psi\|_2^2}\left({{\frac{4\|\psi\|_1^{8/3}}{l^{4/3}}+
4l^{2/3}\|\psi\|_1^{2/3}\|\psi\|_2^2}}\right)^{-1/2}
\]
and therefore
\begin{eqnarray*}
&&\|\nabla\psi\|_2^2\ge{C\|\psi\|_2^6}\left({{\frac{4\|\psi\|_1^{8/3}}{l^{4/3}}+
4l^{2/3}\|\psi\|_1^{2/3}\|\psi\|_2^2}}\right)^{-3/2}\ge
C{\|\psi\|_2^6}\left(\frac{\|\psi\|_1^{4}}{l^{2}}+
l\|\psi\|_1\|\psi\|_2^3\right)^{-1} \\
&&~~~~~~~~
\ge \frac{Cl^2\|\psi\|_2^6}{\|\psi\|_1^{4}+l^3\|\psi\|_1\|\psi\|_2^3}.
\end{eqnarray*}
Hence (\ref{nash-strip}) indeed holds.

We insert (\ref{nash-strip}) into the inequality (\ref{2.9}) and using
the conservation of the $L^1$-norm of $\psi$ (recall that the initial
data is non-negative) obtain
\begin{equation}\label{2.26}
\frac{d\|\psi\|_2}{dt}\le -\frac{C
l^2\|\psi\|_2^5}{\|\psi_0\|_1^4+ l^3\|\psi_0\|_1\|\psi\|_2^3}.
\end{equation}
Integrating (\ref{2.26}) in time we have
\[
C l^2t\le \frac{\|\psi_0\|_1^4}{4\|\psi\|_2^4}+
\frac{l^3\|\psi_0\|_1}{\|\psi\|_2}\le
\frac{1}{z(t)}\left[l^3+\frac{1}{4z^3(t)}\right],
\]
where $z(t)=\|\psi(t)\|_2/\|\psi_0\|_1$, and thus
\begin{equation}\label{z-ineq}
\frac{4z^4(t)}{1+4l^3z^3(t)}\le \frac{1}{Cl^2t}.
\end{equation}
The function on the left side of (\ref{z-ineq}) is monotonically increasing
and hence we have
\begin{equation}\label{l1-l2-n}
\|\psi(t)\|_{2}\le n(t)\|\psi_0\|_1,
\end{equation}
where $n(t)$ is the solution of (\ref{n-eq}).

Let us denote by ${\cal P}_t$ the solution operator for (\ref{v-eq}):
$\psi(t)={\cal P}_t\psi_0$. Then (\ref{l1-l2-n}) implies that
$\|{\cal P}_t\|_{L^1\to L^2}\le n(t)$. The adjoint
operator ${\cal P}_t^*$ is the solution operator for
\begin{eqnarray}\label{v-eq*}
&&\pdr{\tilde\psi}{t}-v\cdot\nabla\tilde\psi=\Delta\tilde\psi\\
&&\tilde\psi(0,x)=\tilde\psi_0(x),~~~x\in{\mathbb R}^d\nonumber
\end{eqnarray}
Note that the preceding estimates rely only on the skew adjointness of
the convection operator $v\cdot\nabla$. Therefore we have the bound
$\|{\cal P}_t^*\|_{L^1\to L^2}\le n(t)$ and hence
$\|{\cal P}_t\|_{L^2\to L^\infty}\le n(t)$ so that
\begin{equation}\label{l2linfty}
\|\psi(t)\|_{L^\infty}\le n(t/2)\|\psi(t/2)\|_{L^2}
\le  n^2(t/2)\|\psi_0\|_{L^1}
\end{equation}
and the proof of Lemma \ref{lemma-l1infty} is complete. $\Box$

A very similar argument leads to an estimate for the solution of
(\ref{small-adv-diff})  with  the $2\pi l$-periodic or
zero  Dirichlet boundary conditions in both $x$ and $y.$ We state
this variant which we will need.

\begin{lemma}\label{periodicdecay11}
Consider equation (\ref{v-eq}) with the $2\pi l$-periodic
boundary conditions in $x$ and $y.$  Assume that the initial data
$\psi_0(\vx)$ is mean zero: $\int \psi_0(\vx)\,d \vx=0.$  Then
there exists a constant $C>0$ such that
\[ \|\psi(t)\|_{L^\infty(D)}\le {C}n^2(t)\|\psi_0\|_{L^1(D)}, \]
where $D=[0,2\pi l]_x\times[0,2\pi l]_y$. Here $n(t)$ is the unique
solution of
\[ \frac{4n^4(t)}{1+4n^3(t)l^3}=\frac{C_1}{ l^2t,} \]
and the constants $C,C_1$ do not depend on $v.$ The same result
holds for the zero Dirichlet boundary conditions, without the
assumption that the initial data is mean zero.
\end{lemma}

\subsection{Variation of temperature on streamlines and gradient
bounds}

The next step of the proof is to estimate the time-space averages
of the oscillations of the solution along the streamlines of $u$
and show that they are small if the advection amplitude is
sufficiently strong. This does not have to be true in general, but
a bound of this type holds if the initial data is uniform on the
streamlines -- then one has only to show that no large
oscillations appear at later times. First, we reduce the problem
to such initial data with the help of the Nash inequality in Lemma
\ref{lemma-l1infty}. The maximum principle implies that it
suffices to prove Theorem \ref{thm2} for the initial data of the
form $T_0=1$ for $|x|\le L_0$ and $T_0=0$ elsewhere. We split the
initial data as
\begin{equation}\label{variation-split}
T_0(\vx)=T_0(\vx)\eta(h(\vx)/\delta_0)+
T_0(\vx)(1-\eta(h(\vx)/\delta_0))=\phi_{01}+\phi_{02}.
\end{equation}
The small parameter $\delta_0$ will be specified later. The cutoff
function $\eta(h)$ satisfies
\[
\hbox{$0\le\eta(h)\le 1$ for all $h\in\Rm$, $\eta(h)=0$ for $|h|\le
1$,
$\eta(h)=1$ for $|h|\ge 2$.}
\]
We split the solution $\phi(t,\vx)$ of (\ref{small-linear}) as a
sum $\phi=\phi_1+\phi_2$. The functions $\phi_{1,2}$ satisfy
(\ref{small-linear}) with the initial data $\phi_{01}$ and
$\phi_{02}$, respectively.
The function $\phi_{01}$ is smooth and is constant (1
or 0)  on the streamlines of the flow. In particular it
is equal to zero in the whole ``water-pipe'' system of
boundary layers around the separatrices. The function
$\phi_{02}$ satisfies a bound
\begin{equation}\label{bd-phi02}
\|\phi_{02}\|_{L^1(D)}\le \left|\left\{\vx\in
D=\IR\times [-\pi l, \pi l]:~ |x|\le L_0,~|h(\vx)|\le
2\delta_0\right\}\right|\le C\delta_0 L_0
\ln(l/\delta_0).
\end{equation}
Therefore, Lemma \ref{lemma-l1infty} implies that if $L_0$ is sufficiently large
and we require that
\begin{equation}\label{delta0}
\delta_0\le \frac{C\theta_0 l^2}{ L_0 (\ln(L_0/l))^2},
\end{equation}
with an appropriate constant $C,$ then the function $\phi_2$
satisfies a uniform upper bound
\begin{equation}\label{var-phi2}
\|\phi_2(t=l^2)\|_{L^\infty(D)}\le \frac{\theta_0'}{10}.
\end{equation}
 Hence
we choose $\delta_0$ as in (\ref{delta0}) and concentrate on the
function $\phi_{1}$ that is initially uniform along the streamlines
of $u$. We drop the subscript one to simplify the notation wherever
this causes no confusion.


Next, we obtain some uniform estimates for solutions
of (\ref{small-linear})
that are initially constant on the streamlines.
\begin{lemma}\label{streamconst}
For any time $t>0$ we have
\begin{equation}\label{gradbound}
\int_0^t \int_D |\nabla \phi|^2 d\vx ds \leq \int_D |\phi_0(x)|^2
d\vx.
\end{equation}
Assume in addition that the initial data $\phi_0(x)$ for the
equation (\ref{small-linear}) are constant on streamlines. Then
\begin{equation}\label{flowinside}
\int_0^t \int_D |u \cdot \nabla \phi|^2 d\vx ds \leq CA^{-2}t \int_D |\Delta \phi_0|^2
d \vx + CA^{-1} \int_D | \phi_0 |^2 d\vx.
\end{equation}
\end{lemma}
{\bf Proof.} Multiplying (\ref{small-linear}) by
$\phi(t,\vx)$ and integrating we trivially obtain
\[ \int_D |\phi(t,\vx)|^2d\vx+\int_0^t\int_D |\nabla\phi(s,\vx)|^2d\vx ds=
\int_D |\phi_0(\vx)|^2d\vx, \]
which implies (\ref{gradbound}).
Next, we multiply (\ref{small-linear}) by $u\cdot\nabla\phi$ and integrate
over $D$ to get
\begin{eqnarray*}
\!\!\int_D |u \cdot \nabla \phi|^2d\vx&= &A^{-1}\int_D u
\cdot
\nabla \phi (\Delta\phi-\phi_t)d\vx\\
&\le&
\frac{1}{2A^2}\int_D \phi_t^2d\vx+\frac 12\int_D
 |u \cdot \nabla
\phi|^2d\vx
-A^{-1} \int_D \nabla \phi^T \nabla u
\nabla \phi d\vx.
\end{eqnarray*}
Thus we get
\begin{equation}\label{small-2}
\int_D \left|u \cdot \nabla \phi \right|^2d\vx\le CA^{-2}\int_D
\phi_t^2d\vx+ C A^{-1} \int_D |\nabla\phi|^2d\vx.
\end{equation}
From (\ref{small-2}) and (\ref{gradbound}) we obtain that
\begin{equation}\label{small-3}
\int_0^t\int_D \left|u \cdot \nabla \phi(\vx, s)\right|^2d\vx
ds\le CA^{-2}\int_0^t\int_D \phi_s^2d\vx ds+ C A^{-1}
\int_D|\phi_0|^2d\vx.
\end{equation}
Recall that initially $\phi_0(x,y)$ is constant on streamlines of
$u$, so that $u\cdot\nabla\phi_0=0$. Since the function $\phi_t$ satisfies the same
equation (\ref{small-linear}) as $\phi$, this implies that
\[
\int_D |\phi_t(t,\vx)|^2d\vx+\int_0^t\int_D |\nabla\phi_t(s,\vx)|^2d\vx ds=
\int_D |\phi_t(0,\vx)|^2d\vx=\int_D |\Delta\phi_0|^2d\vx.
\]
Therefore, bounding from above the time derivative term in
(\ref{small-3}), we arrive at (\ref{flowinside}). $\Box$

Note that the initial data for the function $\phi_1$ in (\ref{variation-split})
obeys an upper bound
\begin{eqnarray*}
&&\int_D|\Delta\phi_{01}|^2d\vx=\int_{-L_0}^{L_0}\int_0^{2\pi
l}\left|\Delta\left[
\eta\left(\frac{h(\vx)}{\delta_0}\right)\right]\right|^2dy dx\\
&&= \int_{-L_0}^{L_0}\int_0^{2\pi l}\left|\frac{\Delta
h(\vx)}{\delta_0}
\eta'\left(\frac{h(\vx)}{\delta_0}\right)+\frac{|\nabla
h(\vx)|^2}{\delta_0^2}
\eta''\left(\frac{h(\vx)}{\delta_0}\right)\right|^2dxdy\le
\frac{CL_0}{\delta_0^3}\ln\left(\frac{l}{\delta_0}\right).
\end{eqnarray*}
Hence, according to (\ref{flowinside}), the total oscillation on streamlines
is bounded by
\begin{eqnarray}\label{in-oscill}
&&\int_0^t\int_D \left|u \cdot \nabla \phi(s,\vx)\right|^2d\vx ds\le
CA^{-2} t\int_D |\Delta\phi_{01}|^2d\vx+ CA^{-1} \int_D
|\phi_{01}|^2d\vx\\
&&~~~~~~~~~~~~~~~~~~~~~~~~~~~~~~~~~~
\le \frac{CA^{-2} \tau L_0}{\delta_0^3}\ln\left(\frac{l}{\delta_0}\right)
+ CA^{-1} L_0\nonumber
\end{eqnarray}
for $t\le \tau$.

\subsection{Boundary layer and cell-to-cell heat conduction}

Lemma~\ref{streamconst} suggests that for most times, there should
be very little temperature variation along the streamlines. Our
next goal is make this statement more precise. Let us fix a time
$\tau$ to be chosen later. We will denote by $h$ and $\theta$ the
coordinates inside each cell, with $h$ being our usual stream
function and $\theta$ the orthogonal coordinate normalized by the
condition $|\nabla \theta|=|\nabla h|$ along the cell boundary and
increasing in the direction of the flow. Let $\delta>0$ be
an
arbitrary,  small number to be chosen later.

\begin{lemma}\label{streamcontrol}
There exists $h_0$ satisfying $2\delta > h_0 > \delta$ such that
\begin{equation}\label{streamvar}
\int_0^\tau \sum\limits_{cells} {\rm sup}_{\theta_1, \theta_2}
|\phi(h_0, \theta_1, t) - \phi(h_0, \theta_2, t)|^2 \,dt \leq
C\delta^{-1}A^{-1}L_0 l\ln\left(\frac{l}{\delta}\right)(A^{-1}\tau
\delta_0^{-3}\ln\left(\frac{l}{\delta_0}\right)+ 1)
\end{equation}
where $\delta_0$ satisfies  (\ref{delta0}).
As a consequence, given a small number $\gamma>0,$ for all times except for a set
$S_\gamma$ of Lebesgue measure at most $\gamma \tau$ we have
\begin{equation}\label{streamvar2}
\sum\limits_{cells} {\rm sup}_{\theta_1, \theta_2} |\phi(h_0,
\theta_1, t) - \phi(h_0, \theta_2, t)|^2 \leq
C\delta^{-1}A^{-1}\gamma^{-1}L_0 l\ln\left(\frac{l}{\delta}\right)(A^{-1}
\delta_0^{-3}\ln\left(\frac{l}{\delta_0}\right)+ \tau^{-1}).
\end{equation}
\end{lemma}
{\bf Proof.}
It is easy to check that
\[ \frac{\partial \phi}{\partial \theta} = |\nabla h|^{-1} |\nabla \theta|^{-1} u \cdot \nabla \phi. \]
Let us denote by $S_a$ the streamline $h=a\in [\delta,2\delta]$ in a
given cell. Then we have
\begin{eqnarray*}
&&{\rm sup}_{\theta_1, \theta_2} |\phi(h,\theta_1) -
\phi(h,\theta_2)|^2 \leq \left( \int_{S_h} |u \cdot \nabla
\phi| \frac{d\theta}{|\nabla \theta||\nabla h|} \right)^2 \\
&&\le C\int_{S_h} |u \cdot \nabla \phi|^2 \frac{d\theta}{|\nabla
\theta||\nabla h|}\int_{S_h} \frac{d\theta}{|\nabla \theta||\nabla
h|} \leq Cl\ln{\frac{l}{\delta}}\int_{S_h} |u \cdot \nabla
\phi|^2 \frac{d\theta}{|\nabla \theta||\nabla h|}.
\end{eqnarray*}
The last step follows from the estimate
\[
\int_{S_h} \frac{d\theta}{|\nabla \theta||\nabla
h|} \leq Cl\ln\left({\frac{l}{\delta}}\right)
\]
 in the tube $\delta < h < 2\delta.$ Integrating
in $h$ and in time and summing over all cells, we obtain
\begin{equation}\label{almostthere1}
\int_0^\tau \int_\delta^{2\delta} \sum\limits_{cells} {\rm
sup}_{\theta_1, \theta_2} |\phi(h, \theta_1, t) - \phi(h,
\theta_2, t)|^2 \,dhdt \leq C l\ln\left(\frac{l}{\delta}\right) \int_0^\tau
\int_{D} |u \cdot \nabla \phi|^2\, d\vx dt.
\end{equation}
Now (\ref{streamvar}) and (\ref{streamvar2}) follow from (\ref{in-oscill})
by an application of the mean value theorem.
$\Box$

We see that for all times but a small set of "exceptional" times,
the value of the temperature on the $h = h_0$ streamline is close
to some constant in any cell. Our next goal is to establish
a control on how different these constants can be for
two neighboring cells. Consider two such cells, $\cC_-$
and
$\cC_+,$ and let $B$ be their common boundary. Let
us choose the coordinates on these two cells so that
$h=0$ on $B,$ $h>0$ on the right cell $\cC_+$, $h<0$
on the left cell $\cC_-,$ and the angular
coordinates $\theta_-$ and
$\theta_+$ are equal to zero in the mid-point of $B.$
Denote by
$D_{\pm}$ the region bounded by curves $\theta_{-,+} = \pm y$ and
streamlines $h = \pm h_0,$ where $y \sim l$ is chosen so that
$D_{\pm}$ is sufficiently far away from the corners of the cells
(so that $|\nabla h|,$ $|\nabla \theta_{\pm}| \geq C$ on
$D_{\pm}$). Also denote by $S_{-,+}$ the pieces of streamlines bounding
$D_\pm.$ Let us define the temperature drop between $\cC_-$ and
$\cC_+$ as follows:
\[ |\phi_+ - \phi_-| = {\rm max}\left\{ 0, {\rm min}_{S_+} \phi -
{\rm max}_{S_-}\phi, {\rm min}_{S_-} \phi - {\rm max}_{S_+}\phi
\right\}. \] Note that if the time is not exceptional,
Lemma~\ref{streamcontrol} implies that maximum and minimum of
$\phi$ along $h=h_0$ streamline in any cell differ by at most
$C\delta^{-1}A^{-1}l \ln
(l/\delta)\gamma^{-1}L_0(A^{-1}\delta_0^{-3}\ln(l/\delta_0)+\tau^{-1}).$
Now we are ready to state
\begin{lemma}\label{diff}
For any $\tau >0,$ we have for
$D=\IR\times [-\pi l, \pi l]$
\begin{equation}\label{diffcontrol}
\int_0^\tau \sum\limits_{cells} |\phi_+ - \phi_-|^2 \,dt
\leq C\delta l^{-1} \int_0^\tau \int_D |\nabla \phi|^2 \,d\vx dt \leq CL_0\delta.
\end{equation}
Therefore, given $\gamma>0,$ for all times with an exception of
a set of Lebesgue measure at most $\gamma \tau,$
we have
\begin{equation}\label{diffcontrol1}
\sum\limits_{cells} |\phi_+ - \phi_-|^2 \leq C \delta L_0 \gamma^{-1}\tau^{-1}.
\end{equation}
\end{lemma}
{\bf Proof.}
The proof is straightforward and has already appeared in \cite{KR}. We sketch it here for the sake of completeness.
Clearly,
\[ \left|
\int_{-h_0}^{h_0} \frac{\partial \phi}{\partial h}(h,\theta)\,dh
\right| \geq |\phi_+ - \phi_-|, \] for any $\theta_{\pm}=\theta.$
Since in the region $D_\pm$ we have $|\nabla h|, |\nabla
\theta_\pm| \sim 1,$ integrating in the curvilinear coordinates we obtain
\[ \int_D |\nabla \phi|^2 d \vx \geq
\int_{-y}^y  \int_{-h_0}^{h_0} \left|\frac{\partial \phi}{\partial h}\right|^2\,d\theta dh \geq
C\delta^{-1}l|\phi_+ - \phi_-|^2, \]
finishing the proof. $\Box$

Lemmas~\ref{streamcontrol} and \ref{diff} allow to
control how much the temperature changes from cell to
cell. Our next lemma summarizes this in a way that will
prove useful.
\begin{lemma}\label{manycells}
Assume that at a certain time $t,$ the estimates (\ref{streamvar2}) and (\ref{diffcontrol1}) hold.
Suppose that there exists a cell ${\mathcal C}$ such that $\phi(h_0, \theta) \geq \beta>0$ at some point
$(h_0,\theta) \in {\mathcal C}$.
Then for at least $N$ cells, we have $\phi(h_0, \theta) \geq \beta/2$ for any $\theta$ inside these cells,
where the number $N$ can be estimated from below as follows:
\begin{equation}\label{manycellsest}
N \geq C\beta^2 \gamma \tau \left( \delta L_0 +
\delta^{-1}A^{-2}L_0 \delta_0^{-3}l\ln(l/\delta) \tau\ln(l/\delta_0) +
\delta^{-1}A^{-1}L_0l\ln(l/\delta) \right)^{-1}.
\end{equation}
\end{lemma}
{\bf Proof}
Consider a cell ${\mathcal C}_1$ which is the closest to ${\mathcal C}$ and
such that there exists a point $(h_0,\theta) \in
{\mathcal C}_1$ with $\phi(h_0,\theta) < \beta/2.$
Then we must have
\begin{equation}\label{droplb}
\sum\limits_{cells} |\phi_+ - \phi_-|
+\sum\limits_{cells} {\rm sup}_{\theta_1,\theta_2} |\phi(h_0, \theta_1)
-\phi(h_0, \theta_2)| > \beta/2,
\end{equation}
where the sum is over all cells between ${\mathcal C}$ and ${\mathcal C}_1.$
On the other hand, Lemmas~\ref{streamcontrol} and
\ref{diff} imply that
\begin{eqnarray}\label{dropcont}
&&\sum\limits_{cells} |\phi_+ - \phi_-|^2 + \sum\limits_{cells} {\rm sup}_{\theta_1,\theta_2} |\phi(h_0, \theta_1)
-\phi(h_0, \theta_2)|^2 \\
&&\leq C\gamma^{-1}\tau^{-1} (\delta L_0 +
\delta^{-1}A^{-1}L_0l\ln(l/\delta)(A^{-1}\delta_0^{-3}\tau\ln(l/\delta_0)
+1)).\nonumber
\end{eqnarray}
Since $N \sum_{m=1}^N a_n^2 \geq \left( \sum_{m=1}^M a_n \right)^2,$
a combination of (\ref{droplb}) and (\ref{dropcont}) yields
(\ref{manycellsest}).
$\Box$

Now we can explain the strategy of the proof of Theorem~\ref{thm2}
in more detail. Assume that at a certain "good" (not exceptional
in the sense of Lemmas~\ref{streamcontrol}, \ref{diff}) time we
have a sufficiently high temperature $\phi = \beta$ in some cell
on a streamline $h = h_0.$ Then an appropriate choice of $\delta$
(hence $h_0$) and an application of Lemma~\ref{manycells} ensure
that $\phi \geq \beta/2$ in many cells. Assume that for a
sufficiently large portion of times $\leq \tau$, the temperature
at $h=h_0$ streamlines is high in many cells. Clearly, the
interiors of the cells will heat up too, giving (for a suitable
choice of $\delta$ ensuring $N$ is large enough) a contradiction
with the preservation of the $L^1$ norm of $\phi.$ On the other
hand, if the temperature on the streamline $h=h_0$ is low most of
the time, we expect the solution inside the cells to be small and
quenching to happen if the cells are sufficiently small -- the
last condition ensures that the interaction between the boundary
and the interior of the cell happens on a short time scale. To
make the above plan work, we need a good understanding of an
auxiliary "cell heating" problem. The next section is devoted to
this goal.

\subsection{An auxiliary cell heating problem}

Consider a domain $\Omega$ which is a sub-domain of a cell, defined by
 the condition $|h| \geq |h_0|$, that is, $\Omega$ is
the {\em  interior}
of the closed streamline $\{h=h_0\}$. We consider the following
initial-boundary value problem on $\Omega:$
\begin{eqnarray}\label{auxcell}
w_t +Au \cdot \nabla w - \Delta w =0 \\
w(t,\vx) = \sigma (t), \,\,\,\vx \in \partial \Omega \nonumber \\
w(\vx,0) = g(\vx). \nonumber
\end{eqnarray}
Here $\sigma (t)$ is some (smooth) function of time, independent
of $\vx.$ Let us derive a compact formula for the solution of
(\ref{auxcell}). Set $v = w - \sigma(t).$ Then $v$ satisfies the
zero Dirichlet boundary conditions, and we have
\[ v_t +A u \cdot \nabla v - \Delta v = -\sigma'(t),
 \,\,\,\, v(\vx,0)=g(\vx)-\sigma(0). \]
Let us denote by $\cH$ the differential operator $-\Delta + A u \cdot
\nabla$ with the zero Dirichlet boundary conditions on $\partial \Omega$
and by $e^{-t\cH}$ the positive semigroup generated by $\cH.$
Using the Duhamel formula, we get
\begin{equation}\label{auxsol22}
v(t,\vx) = e^{-t\cH}(g-\sigma(0){\bf 1}) - \int_0^t e^{-(t-s)\cH}{\bf 1}
\sigma'(s)\,ds.
\end{equation}
Here the semigroup is applied to a function equal identically to one
in all of $\Omega,$ denoted by ${\bf 1}.$ Integrating by parts in
(\ref{auxsol22}), and recalling that $w=v+\sigma(t),$ we obtain
\begin{equation}\label{cellsol1}
w(t,\vx) = e^{-t\cH}g(\vx) + \int_0^t \sigma(s) \cH
e^{-(t-s)\cH} {\bf 1} \, ds.
\end{equation}

It will be convenient to use probabilistic interpretation of
(\ref{cellsol1}). Recall that we can represent the semigroup
action in the following way (see, e.g. \cite{Freidlin1}):
\[ e^{-t\cH}f(\vx) = \IE^\vx \left[ f(X_t^\vx(\omega))\right],
\] where $X^\vx_t$ is the diffusion process corresponding to $\cH$ starting at
$\vx,$ and the expectation $\IE^\vx$ is taken over all paths that
start at $\vx$ and never leave $\Omega$ up to time $t.$
The latter restriction corresponds to the zero Dirichlet boundary
condition on $\partial \Omega.$ Note that in the case $f \equiv
1,$ the function $Q(t,\vx) = e^{-t\cH}{\bf 1}$ coincides with
probability that the stochastic process $X^\vx_t(\omega)$ never
leaves $\Omega$ before time $t:$
\[ e^{-t\cH}{\bf 1} = \IE^\vx\left[ {\bf 1} (X^\vx_t(\omega))\right]
= \IP(X^\vx_s(\omega) \in \Omega, \,\,\, \forall \,\,s \leq t). \]
In terms of the function $Q(t,\vx)$ expression (\ref{cellsol1}) becomes
\begin{equation}\label{cellsol2}
w(t,\vx) = e^{-t\cH}g(\vx) - \int_0^t \sigma(s)
\partial_t Q(t-s,\vx) \, ds
\end{equation}
In order to control the solution of an auxiliary cell problem, we need
to control the properties of $Q(t,\vx)$, uniformly in the flow
amplitude $A$ -- this is akin to Lemma \ref{lemma-l1infty}.  The first
result we need in this direction is the following lemma, showing that
the exit probability is bounded from above by a constant independent
of $A.$
\begin{lemma}\label{lowercell}
For any $t$ satisfying $0 \leq t \leq C_1 l^2$ we have
\begin{equation}\label{lbQ}
\int_\Omega Q(t,\vx) \,dx \geq C_2 l^2 >0,
\end{equation}
where $C_2$ is a universal constant depending only on $C_1$ but
not on $A$ or $l.$
\end{lemma}
{\bf Proof.} Let us rewrite the generator $\cH$ in the natural
coordinates $h,\theta:$
\[ \cH = |\nabla \theta|^2 \frac{\partial^2}{\partial \theta^2}
+  |\nabla h|^2 \frac{\partial^2}{\partial h^2}+ \Delta h
\frac{\partial}{\partial h} + (\Delta \theta - A |\nabla \theta|
|\nabla h|) \frac{\partial}{\partial \theta}.
\]
Note that $\Delta h = -2l^{-2}h$ -- this property is by no means
crucial for our analysis but it does simplify some computations. The
diffusion process $X^\vx_t$ corresponding to $\cH$, written in the
$(h,\theta)$-coordinates, is given by
\begin{eqnarray}\label{OUp}
dX^h_t = \sqrt{2} |\nabla h| dB^{(1)}_t - 2 l^{-2} X^h_t dt \\
dX^\theta_t = \sqrt{2} |\nabla \theta| dB^{(2)}_t +(\Delta \theta
- A |\nabla h| |\nabla \theta|)dt, \nonumber
\end{eqnarray}
where the values of all functions are taken at a point $(X^h_t,
X^\theta_t)$, while $B^{(1)}$ and $B^{(2)}_t$ are independent one-dimensional
Brownian motions. Clearly, $Q(t,\vx) \geq P(X^{h(\vx)}_s
\in [0, l],\,\,\,\forall \,\,s \leq t).$ It is not
difficult to see that the exit probabilities of $X^h_t$
are majorized by the exit probabilities of the
Ornstein-Uhlenbeck process where the factor $|\nabla
h|$ in (\ref{OUp}) is dropped.  Indeed, let us introduce
\[
\alpha (t,\omega) =
\sqrt{2} |\nabla h (X^h_t(\omega), X^\theta_t(\omega))|.
\]
Multiplying (\ref{OUp}) with $e^{2l^{-2} t}$ and integrating leads
to
\[
e^{2l^{-2}t} X^{h(\vx)}_t =
h(\vx) +  \int_0^t e^{2 l^{-2} s} \alpha(s,\omega) dB^{(1)}_s.
\]
Making a random time change (see, e.g., \cite{Oks}), we find that the
integral above has the same distribution as the Brownian motion
$B^0_{\rho(t,\omega)},$ where
\[ \rho (t,\omega) =  \int_0^t e^{4l^{-2} s} |\alpha (s,\omega)|^2\,ds. \]
Since $|\alpha (s, \omega)|^2 \leq 2,$ we have
\[ \rho(t, \omega) \leq \rho(t) \equiv \frac12 l^2 (e^{4 l^{-2} t}-1), \]
and therefore
\begin{equation}\label{bmbound}
 Q(t,\vx) \geq P( B^{h(\vx)}_{\rho(s)} e^{-2l^{-2} s}
\in [0,l],
\,\,\,\forall \,\,s \leq t).
\end{equation}
Let us remark that the expression on the right hand side is exactly
the exit probability for the Ornstein-Uhlenbeck process mentioned
above. The claim of the lemma now follows from a simple rescaling $\tau=t/l^2$,
$Y_\tau=l^{-1}X_{l^2\tau}$. $\Box$

Next, we need some information on the behavior of $\partial_t
Q(t,\vx).$
\begin{lemma}\label{ubforderiv}
For any $t$ satisfying $C_1 l^2 \geq t \geq 0,$ we have
\[
\int_\Omega \partial_t Q(t,\vx) d\vx \leq -C_3,
\]
with a constant $C_3>0$ that depends on $C_1$ but is independent
of $A$ or $l$. Moreover, $\int_\Omega
\partial_t Q(t,\vx) d\vx$ is monotonically increasing in
time.
\end{lemma}
{\bf Proof.}
According to the previous lemma,
\[
\int_\Omega Q(t,\vx) \,d\vx \equiv \|Q(t,\vx)\|_{L^1}
\geq C_2l^2 >0
\]
for any $t \leq C_1 l^2.$ The Cauchy-Schwartz inequality then implies a lower
bound for the $L^2$ norm: $\|Q(t,\vx)\|_{L^2} \geq C_2
l.$ Recall that the function $Q(t,\vx)$ solves
\[
\pdr{Q}{t}= \cH Q(t,\vx)
\]
with the
zero  Dirichlet boundary conditions and the initial condition
$Q(0,\vx)=1.$ Then
\begin{equation}\label{l2decay}
\frac{1}{2}\partial_t \|Q\|^2_{L^2} = - \|\nabla Q\|^2_{L^2} \leq
-l^{-2} \|Q\|^2_{L^2}
\end{equation}
by the Poincare inequality. Thus we have $\int Q Q_t d\vx \leq
-C_2^2$ where $C_2$ is independent of
$A$ and $l$. It is clear
that
$0
\leq Q(t,\vx)
\leq 1,$ and also that
$\partial_t Q(t,\vx) \leq 0$ since $Q(t,\vx)$ is the
probability that the diffusion process starting at
$\vx$ does not exit $\Omega$ before time $t.$ This
implies the first statement of the lemma.

To prove the second statement, we use again the monotonic decay of
$Q(t,\vx)$ in time.  Fixing some time $t>0$ and integrating by parts
(using the fact that the flow $u$ is tangent to $\partial \Omega$) we
find that
\[
\int_\Omega \partial_t Q(t,\vx) \,d\vx =
\int_{\partial \Omega} \frac{\partial Q}{\partial n}\, ds.
\]
Since the function $Q$ vanishes on the boundary and decays (in
time) inside, we see that $\int_\Omega \partial_t Q \, d\vx$ is
monotonically increasing in time. $\Box$

\begin{corollary}\label{l1lb}
The solution $w(t,\vx)$ of (\ref{auxcell}) with any
nonnegative boundary and initial data satisfies
\begin{equation}\label{eql1lb}
\int_\Omega w(t,\vx)d\vx \geq C \int_0^t \sigma (s)\,ds,
\end{equation}
where $C$ is a positive constant depending on $t$ but not on $A.$
If $t \leq C_1l^2,$ the constant $C$ can be chosen
independent of
$l.$
\end{corollary}
{\bf Proof.} The corollary follows immediately from the previous lemma and the
representation (\ref{cellsol2}). $\Box$

Since we have no control over what happens at a small set of
exceptional times, we need an estimate from above on how much the
$L^1$ norm of the solution can change if the boundary data is
close to one for a short time. The second part of
Lemma~\ref{ubforderiv} implies that it is sufficient to look at
the situation when this hot period occurs at the end of the
interval $[0,t]$. More precisely, if we replace $\sigma(t)$ by its
monotonically increasing rearrangement then the $L^1$-norm of the
solution increases. The following lemma will be useful in such
scenario.
\begin{lemma}\label{extimin}
For a time $t$ satisfying $0 < t < l^2,$ we have
\begin{equation}\label{smalltimedrop}
\int_\Omega (1-Q(t,\vx))\,d\vx \leq C l^2 (t
l^{-2})^{1/2}
\ln
\left(\frac{l^2}{t}\right).
\end{equation}
\end{lemma}
{\bf Proof.}
Note that
\[ \int_\Omega (1-Q(t,\vx))\,d\vx =
\int_\Omega \IP(\exists \,\,s< t: \,\,\,X^\vx_s(\omega) \notin
\Omega)\, d \vx.
\]
It follows from the proof of Lemma~\ref{lowercell}
that the probability on the right hand side is majorized by
\begin{equation}\label{probo}
\IP( \exists \,\,s< t: B^{h(\vx)}_{\rho(s)} e^{-2 l^{-2} s} \notin
[0,l]) \leq \IP( \exists \,\,s< t: B^{h(\vx)}_{\rho(s)} \notin [0,l]),
\end{equation}
where $\rho(t)=\frac12 l^2 (e^{4l^{-2} t}-1).$ Take a constant $C$ so
that $\rho(t) \leq C t$ for any $t$ as in the statement of the
lemma. Let $\tilde{t}= t l^{-2}$ be the rescaled time. Then the
probability on the right hand side of (\ref{probo}) does not exceed
the probability that $B^0_{s}$ leaves the interval $C^{-1/2}[-h/l, 1-h/l]$
before the time $\tilde{t}.$ Using the reflection principle, we can estimate
this probability from above by one, if $h/l$ or $1-h/l$ are less than
$\sqrt{\tilde{t}},$ and by $C_n (\tilde{t}l/h)^n$ otherwise (where $n
>0$ is arbitrary). Integrating and taking into account the Jacobian of
transformation between $h, \theta$ and $\vx,$ we arrive at
\[
\int_\Omega (1-Q(t,\vx))\,d\vx\le Cl(\tilde
tl^{2})^{1/2}
\ln\left(\frac{Cl}{(\tilde tl^{2})^{1/2}}\right)=Clt^{1/2}
\ln\left(\frac{Cl^2}{t^{1/2}}\right),
\]
which is nothing but (\ref{smalltimedrop}). $\Box$

\subsection{Proof of Theorem~\ref{thm2}}

Now we are ready to complete the proof of Theorem~\ref{thm2}. In
this proof, we will work with the time scale $\tau
=C_0l^2$ for a sufficiently large universal constant
$C_0.$ It means that in all estimates of the previous
sections that we are going to use, we take $\tau =
C_0l^2.$ This ensures, in particular, applicability of
Lemma~\ref{ubforderiv} and Corollary~\ref{l1lb}. If we
show that $\phi(t,\vx) \leq \theta_0'/2=
\theta_0/(2e)$ for some $t
\leq
\tau,$ this will be
sufficient for quenching since we assume that
$\tau \leq \tau_c = M^{-1}$ (small cell size).

{\bf Proof.} We have a set $S_\gamma$ of exceptional times of size
at most $\gamma \tau$ such that on the complement  $S_\gamma^c$ of
$S_\gamma$ the estimates (\ref{streamvar2}), (\ref{diffcontrol1})
are valid. The value of the constant $\gamma$ will be chosen later
and will be independent of $A$ and $l.$ Assume first that for all
$t \in S_\gamma^c$ except for a set $S_b$ (of ``bad'' times)  of
size at most $\gamma \tau$ we have $\phi(h_0, \theta) \leq \gamma$
for all cells and all $\theta.$ We are going to show that if
$\gamma$ is chosen sufficiently small (with the choice being
uniform in $A,l$), then we must have quenching. Let $\Omega_n$ be
the subset of a cell ${\mathcal C}_n$ enclosed by $|h| = |h_0|.$
Then for any $\vx \in \Omega_n$ we have $\phi(t,\vx) \leq
\tilde{\phi}(t,\vx)$ where $\tilde{\phi}(t,\vx)$  satisfies
\[
\tilde\phi_t-\Delta \tilde{\phi} + A u \cdot \nabla \tilde{\phi} =0
\]
with the initial data
$\tilde{\phi}(\vx,0) = 1$ in $\Omega_n$,
and the boundary data given by
\[
\tilde{\phi}(s, \vx) |_{h=h_0} =\left\{
\begin{matrix} 1, & s \in S_\gamma \cup S_b; \cr
\gamma, & {\rm otherwise}
\cr\end{matrix}\right. .
\]
This bound on $\phi$ follows from the representation formula
(\ref{cellsol2}) and the above assumption on the behavior of
$\phi(h_0,\theta).$ By linearity, inside each region $\Omega_n,$
we have $\tilde{\phi}(t,\vx) =
\tilde{\phi}_1(t,\vx)+\tilde{\phi}_2(t,\vx),$ where
$\tilde{\phi}_1$ satisfies the zero Dirichlet boundary conditions,
while $\tilde{\phi}_2$ has zero initial data. Let $t_1 = C_1 l^2,$
where $C_1$ is a universal constant which will be chosen below. A
simple argument using (\ref{l2decay}) shows that
\begin{equation}\label{l1dec}
 l^{-2} \|\tilde{\phi}_1
(t_1,\cdot)\|_{L^{1}(\Omega_n)} \leq
l^{-1}\|\tilde\Phi_1(t_1,\cdot)\|_{L^2(\Omega_n)}\rightarrow
0
\end{equation}
uniformly in $A,l$ as $C_1 \to \infty$.

The function $\tilde{\phi}_2$ can be estimated in the $L^1$-norm
using (cf. (\ref{cellsol2}))
\begin{eqnarray}\nonumber
&&\tilde\phi_2(t,\bx)=-\int^t_0\sigma(s)\partial_t
Q(t-s,\vx)ds
\end{eqnarray}
and the second statement in
Lemma~\ref{ubforderiv} as
\begin{eqnarray}\nonumber
&&\|\tilde \phi_2(t_1,\cdot)\|_{L^1(\Omega_n)}
\leq -\int^{t_1}_0\sigma(s)
\int_{\Omega_n}\partial_{t_1} Q(t_1-s,\vx)d\vx ds\\
&&~~~~~~~~~~~~~\leq-\int_{t_1-|S_\gamma\cup
S_b|}^{t_1}\int_{\Omega_n}\partial_{t_1} Q(t_1-s,\vx)d\vx
ds-\gamma\int^{t_1-|S_\gamma\cup S_b|}_{0}\int_{\Omega_n}
\partial_t Q(t_1-s,\vx)d\vx ds\nonumber\\
&&~~~~~~~~~~~~\leq \int_{\Omega_n} \left[1 - Q(|S_\gamma\cup S_b|,
\vx)\right]d\vx + \gamma \int_{\Omega_n} \left[1-
Q(t_1,\vx)\right]d\vx \nonumber\\
&&~~~~~~~~~~~~\leq Cl^2 \left[ (C_1 \gamma)^{1/2} \ln
\left(\frac{1}{C_1\gamma}\right) + C_1 \gamma\right]
\label{boundcon}
\end{eqnarray}
by Lemma~\ref{extimin}.

It follows now from the parabolic maximum principle  that for all times $t \geq
t_1,$ the solution $\phi(t,\vx)$ of the original linear advection-diffusion
problem (\ref{small-linear}) in any cell ${\mathcal C}_n$ satisfies
\[
\phi(t,\vx) \leq \phi^*(t,\vx)+ \phi_{2}(t,\vx).
\]
The function  $\phi_2$ has been
estimated in   (\ref{bd-phi02})
and (\ref{var-phi2}). The function $\phi^*(t,\vx)$ solves
(\ref{auxcell})  for $t\geq t_1$ with the $2\pi
l$-periodic boundary conditions in
$x,y$ and at time $t=t_1=C_1l^2$ is given by
\[
\phi^*( t_1,\vx) = \left\{ \begin{array}{ll}
\tilde{\phi}(t_1,\vx), & \vx \in \Omega_n \\
1 & \vx \in {\mathcal C}_n \setminus \Omega_n. \end{array} \right.
\]
Lemma~\ref{periodicdecay11} allows us to choose a positive
number $\tilde{\gamma}$ so that if
\[
\|\phi^*(t_1,\vx) \|_{L^1({\mathcal C}_n)} \leq
\tilde{\gamma}l^2,
\]
then
\[
\|\phi^*( t_1+l^2,\vx) \|_{L^\infty({\mathcal C}_n)}
\leq Cn^2(l^2) \|
\phi^*(t_1,\vx) \|_{L^1({\mathcal C}_n)} \leq
\theta_0'/4.
\]
This is possible, as $n^2(l^2)\sim l^{-2}$ (see Remark after
Lemma \ref{lemma-l1infty}).
Note that
\[ \| \phi^*( t_1,\vx) \|_{L^1({\mathcal C}_n)} \leq
\|\tilde{\phi}_1(t_1,\vx) \|_{L^1(\Omega_n)} +
\|\tilde{\phi}_2(t_1,\vx) \|_{L^1(\Omega_n)} + C\delta l
\ln(l/\delta). \] Choosing $C_1$ large enough we can use
(\ref{l1dec}) to estimate the first term on the right hand side
and make sure it does not exceed $\tilde{\gamma}l^2/3$. Next we
choose $\gamma$ small enough so that (\ref{boundcon}) gives
similar control of the second term. Finally, choose $\delta$ small
enough so that the last term is also sufficiently small. With this
choice of $C_1,$ $\gamma$ and $\delta$ (uniform in $A,l$) we have
\[ \|\phi( t_1+l^2,\vx)\|_{L^\infty} \leq
\| \phi^*(t_1+l^2,\vx) \|_{L^\infty({\mathcal
C}_n)}+\theta_0'/10 \leq
\theta_0'/2
\] and thus quenching.

It remains to consider the case when there exists a set of bad
times $S_b \in [0,t_1]$ of size at least $\gamma \tau$ (with
$\gamma$ universal constant determined above) such that for any $t
\in S_b,$ there exists a cell ${\mathcal C}_n$ such that for some
point $(h_0,\theta) \in {\mathcal C}_n,$ we have $\phi(h_0,\theta)
\geq \gamma.$ We are going to show that this cannot be true if $A$
is large enough, thus forcing the scenario which is  considered
above and leads to quenching. Indeed, Lemma~\ref{manycells}
implies that there are at least
\begin{equation}\label{numcell22}
 N=C\gamma^3 \tau \left( \delta L_0 + \delta^{-1}A^{-2}L_0
 \delta_0^{-3}\tau l\ln(l/\delta)\ln\left(l/\delta_0\right)+
\delta^{-1}A^{-1}L_0 l\ln(l/\delta)\right)^{-1}
\end{equation}
cells such that $\phi(h_0,\theta) > \gamma/2$ for any $\theta$ in
these cells. For each cell, let $\sigma (t) = {\rm min}_\theta
\phi(h_0, \theta, t).$ On each cell, solve the initial-boundary
value problem (\ref{auxcell}) with $g=0,$ and denote the solution
$\overline{\phi}(t,\vx).$ By the parabolic maximum principle, we
have that $\phi(t,\vx) \geq \bar{\phi}(t,\vx)$  for $|h |\geq |
h_0|.$ Applying Corollary~\ref{l1lb}, we obtain
\begin{equation}\label{L1contr}
 \int_D \phi( \tau,\vx)\,d\vx \geq C \sum\limits_{cells}
\int_0^\tau \sigma(s)\,ds = C \int_0^\tau
\sum\limits_{cells} \sigma(s) \,ds \geq C \gamma^2
\tau N,
\end{equation}
where $N$ is given by (\ref{numcell22}). We claim that if $
\delta=\delta
(L_0)$ is chosen to be sufficiently small and $A=A(L_0)$
sufficiently large then (\ref{L1contr}) leads to
\begin{equation}\label{L1normtr}
\int_D \phi(\tau,\vx) \, d\vx \gg l L_0,
\end{equation}
obtaining a contradiction since the $L^1$ norm of $\phi$ is
preserved by evolution. Indeed, since $\gamma$ and $\tau/l^2$ are
fixed constants, it suffices to choose $\delta$ so that
\begin{equation}\label{deltacon}
 \delta \ll \frac{\tau^2}{L_0^2 l} ,
\end{equation}
and then choose $A$ so that
\begin{equation}\label{Acon}
 A \gg{\rm max} \left\{ \frac{C(l)L_0}{ \delta_0^{3/2} \delta^{1/2}}
 \left(\ln\frac{l}{\delta_0}\ln\frac{l}{\delta}\right)^{1/2},
 \frac{C(l)L_0^2 }{ \delta}\ln\left(\frac{l}{\delta}\right) \right\}.
\end{equation}
Recalling the formula (\ref{delta0}) for $\delta_0$, we discover
that $A> C(l)L_0^4\ln L_0$ is sufficient to satisfy (\ref{Acon}),
(\ref{deltacon}), completing the proof of Theorem~\ref{thm2}.
$\Box$

\end{document}